\documentclass[12pt]{amsart}

\input{sommers_macros.tex}
 \usepackage{palatino}

\newcommand{\g}{{\mathfrak g}}
\newcommand{\tor}{{\mathfrak t}}
\newcommand{\bo}{{\mathfrak b}}

\def\Ad{\textup{Ad } }
\def\H{{\mathcal B}}
\def\I{{\mathcal I}}

\def\A{\mathcal A}
\def\hyp{\mathcal A}

\newcommand{\al}{{\alpha}}

\newcommand{\ro}{{\Phi}}
\newcommand{\co}{{\Phi^{\scriptscriptstyle\vee}}}

\newcommand{\posroots}{\ro^+}
\newcommand{\ideal}{\mathcal I}
\newcommand{\nilrad}{\mathfrak{n}}

\newcommand{\hgt}{\text{ ht}}

\newcommand{\colat}{{Q^{\scriptscriptstyle\vee}}}

\begin{document}

\setlength{\topmargin}{.05in}
\setlength{\textheight}{8.6in} 
\setlength{\oddsidemargin}{-.1in}
\setlength{\evensidemargin}{-.1in} \setlength{\marginparsep}{0in}
\setlength{\marginparwidth}{0in}

\title {Exponents for $B$-stable ideals}

\author{Eric Sommers}
\author{Julianna Tymoczko}
\thanks{E.S. was supported in part by NSF grants DMS-0201826 and DMS-9729992}
\thanks{2000 MSC: 20G05 (14M15, 05E15)}

\address{
Department of Mathematics \& Statistics\\
University of Massachusetts--Amherst\\
Amherst, MA 01003 
}

\email{esommers@math.umass.edu}

\address{
Department of Mathematics \\
University of Michigan\\
Ann Arbor, MI  48109-1109
}

\email{tymoczko@umich.edu}

\date{1/06/04; 6/02/04}

\begin{abstract}

Let $G$ be a simple algebraic group over the complex numbers
containing a Borel subgroup $B$.
Given a $B$-stable ideal $I$ in the nilradical of the Lie algebra of $B$,
we define natural numbers $m_1, m_2, \dots, m_k$ which we call
ideal exponents.  We then propose two conjectures where these
exponents arise, proving these conjectures in types $A_n, B_n, C_n$ and some other types.  

When $I = 0$, we recover the usual exponents of $G$ by Kostant \cite{kostant:princ}
and one of our conjectures reduces to a 
well-known factorization of the Poincar{\' e} polynomial of the Weyl group
(see also \cite{macdonald:poincare}).
The other conjecture reduces to a well-known result of Arnold-Brieskorn
on the factorization of the 
characteristic polynomial of the corresponding
Coxeter hyperplane arrangement (see \cite{orlik-terao}).

\end{abstract}
\maketitle

\section{Introduction}

Let $G$ be a simple algebraic group over the complex numbers
containing a Borel subgroup $B$.
The ideals in the nilradical of the Lie algebra of $B$, which are stable 
under the action of $B$, have attracted much recent attention.  

In this paper we define a sequence of natural numbers 
$m_1, m_2, \dots, m_k$ for each $B$-stable ideal $I$, and call them ideal exponents.
The definition is a generalization of the usual exponents of $G$ in the case
where $I = 0$, via Kostant's proof relating the heights of positive roots to the 
exponents \cite{kostant:princ}.

We then conjecture (and prove in type $A_n, B_n, C_n$ and in some other cases) 
two results about these ideal exponents.  
The first concerns a Poincar{\' e} polynomial defined for each ideal which generalizes
the Poincar{\'e} polynomial for the Weyl group.
The conjecture is that this new 
polynomial factors according to the ideal exponents 
just as the usual polynomial factors according to the usual exponents.
This result is relevant for the study of regular nilpotent Hessenberg varieties
(there is one for each ideal) since the combinatorially-defined Poincar{\' e} 
polynomials in this paper should be the actual 
Poincar{\' e} polynomials for these varieties.  This is known in many cases, as studied in \cite{T}.

The second occurrence of these new exponents is in the context of 
a hyperplane arrangement defined for each ideal.  The hyperplane arrangement
in question consists of those hyperplanes which correspond to the positive roots
whose root space does not belong to the ideal.
Generalizing the known result that the usual exponents are the roots of the
characteristic polynomial for the full Coxeter arrangement,
we conjecture (and prove in 
the classical types)
that the characteristic polynomial of this new hyperplane arrangement 
has (non-trivial) roots $m_1, m_2, \dots, m_k$.
We also speculate that these arrangements are free
(which we also prove in the classical types).

The paper concludes with speculation linking these two occurrences of the
ideal exponents.

\section{Notation}

Fix a  maximal torus $T$ in $B$ and let
$(X, \ro, Y, \co)$ be the root datum determined
by $G$ and $T$, and let $W$ be the Weyl group.
Let $\Pi \subset \posroots$ denote the simple roots and positive roots 
determined by $B$.  As usual, $\langle \ , \ \rangle$ denotes the pairing
of $X$ and $Y$.  Let $\colat$ denote the lattice
in $Y$ generated by $\co$ (the coroot lattice). 
We denote the
standard partial order on $\ro$ by $\preceq$; so $\al \preceq \beta$ for 
$\al, \beta \in \ro$ if and only if $\beta - \al$ is a sum of positive roots.
As is customary, we write $\al \prec \beta$ if $\al \preceq \beta$ and $\al \neq \beta$.
For $\beta \in \ro$, write $\beta = \sum_{\al \in \Pi} c_{\al} \al$ and
let $\hgt(\beta) = \sum_{\al \in \Pi} c_{\al}$ denote the height of $\beta$.

We define an ideal (also called an upper order ideal) 
$\ideal$ of $\posroots$ to be a collection
of roots such that if 
$\alpha \in \ideal, \beta \in \posroots$, and $\alpha+\beta \in \posroots$,
then $\alpha + \beta \in \ideal$.
In other words, if $\al \in \ideal$ and $\gamma \in \posroots$ with $\al \preceq \gamma$,
then $\gamma \in \ideal$. 

Let $\g, \bo, \tor$ be the Lie algebras of $G$, $B$, $T$, respectively.
It is easy to see that $B$-stable ideals in the nilradical $\nilrad$ of $\bo$
are naturally in bijection with the ideals of $\posroots$.
Namely, if $I$ is a $B$-stable ideal of $\nilrad$, it is stable under the action of $T$,
hence $I$ is a sum of root spaces.  Denote by $\ideal$ the set of roots 
whose root space is contained in $I$.  Then $\ideal$ is an ideal of $\posroots$ 
and this map is a bijection.      

\section{Ideal exponents}

In this section, motivated by Kostant's proof relating the heights of the positive
roots and the usual exponents of $G$, we define exponents for each ideal.
Our definition is an easy modification:  we consider only those positive roots which
do not lie in the ideal.

For an ideal $\ideal \subset \posroots$, let
$\ideal^c = \posroots - \ideal$ be the positive roots not in $\ideal$.
Define 
$$\lambda_i = \# \{ \al \in \ideal^c \, | \, \hgt( \al) = i \}.$$
We first observe

\begin{prop} \label{partition}
The $\lambda_i$ give a partition of the number of roots in $\ideal^c$. 
That is, $$\lambda_1 \geq \lambda_2 \geq \dots $$
In addition, $\lambda_1 > \lambda_2$.
\end{prop}

\begin{proof}
This is easy to check in the classical groups and was checked on a 
computer in the exceptional
groups.  
\end{proof}

Let $k = \lambda_1$, which is just the number of simple roots in $\ideal^c$.
We define $m^{\ideal}_k \geq  \dots \geq m^{\ideal}_1$ to be the dual partition
of $\lambda_i$.  In other words,
$m^{\ideal}_i = \# \{ \lambda_j \, | \, \lambda_j \geq k - i + 1 \}$.

\begin{defn}
The {\bf ideal exponents} of $\ideal$,
also called $\ideal$-exponents,
are the natural numbers $$m^{\ideal}_k \geq m^{\ideal}_{k-1} \geq \dots \geq m^{\ideal}_1.$$
\end{defn}

It follows from the fact that $\lambda_1 > \lambda_2$ that $m^{\ideal}_1 = 1$.
We also observe, as mentioned previously, 
that when $\ideal = \emptyset$ these are the usual exponents 
(in this case $k$ equals the rank of $G$) \cite{kostant:princ}.  

We suspect that there are many situations where these
new exponents will arise.  We propose two situations in what follows, namely
Theorem \ref{thm1} and Theorem \ref{thm2}.

\section{Poincar{\' e} polynomials for ideals}

Let $R \subset \posroots$ be any subset of the positive roots.  
Given $S \subset R$ we say that 
$S$ is {\bf R-closed} if  
$\al, \beta \in S$ and $\al + \beta \in R$, 
then $\al + \beta \in S$.

Given an ideal $\ideal \subset \posroots$, 
we are interested in those 
subsets $S$  of $\ideal^c$ with the property 
that both $S$ and its complement $S^c:= \ideal^c - S$
are $\ideal^c$-closed.  

These subsets are analogous to Weyl group elements.  Indeed,
if $\ideal = \emptyset$ and if $w \in W$, then 
$$N(w) := \{ \al \in \posroots \, | \, w(\al) \prec 0 \}$$
is $\ideal^c$-closed (in this case, $\ideal^c = \posroots$). 
Conversely every subset of $\posroots$ which is $\posroots$-closed
is equal to $N(w)$ for a unique $w \in W$.  
This is well-known and goes back to \cite{kostant:borel}.

Given this background, we
define a subset $S$ of $\ideal^c$ to be of {\bf Weyl-type} for $\ideal$ 
if both $S$ and $S^c$ are $\ideal^c$-closed.
Let ${\mathcal W}^{\ideal}$ denote the subsets of $\ideal^c$ of Weyl-type.
One of the main results of this paper can now be formulated.

\begin{thm} \label{thm1}
Let $\ideal$ be an ideal in $\posroots$.  Then
in types $A_n, B_n, C_n, G_2, F_4, E_6$
\begin{equation} \label{poincare}
\sum_{S \in {\mathcal W}^{\ideal}} t^{|S|} 
= \prod_{i=1}^{k} (1 +t + t^2 + \dots + t^{m^{\ideal}_i}),
\end{equation}
where the $m^{\ideal}_{i}$ are the exponents of $\ideal$.
\end{thm}

We conjecture that the theorem also holds in the remaining cases.  We defer the proof,
which is case-by-case, until Section \ref{main_proof}.

In the case where $\ideal^c = \posroots$, the theorem is well-known
\cite{kostant:princ}, \cite{macdonald:poincare}.  
On the one hand,
the $\ideal$-exponents become the usual exponents as mentioned previously.
On the other hand, if $l(w)$ denotes the length of $w \in W$ with respect to
the set of simple reflections coming from $\Pi$, then $l(w) = |N(w)|$ and so 
the left-hand side of Equation (\ref{poincare}) is equal to 
$$ \sum_{w \in W} t^{l(w)}.$$
Then Equation (\ref{poincare}) becomes the well-known factorization of the Poincar{\' e} polynomial of 
the Weyl group, which is also the Poincar{\' e} polynomial of the flag variety $G/B$
when $t$ is replaced by $t^2$.

\section{Some old results for general root systems}

Many of the results in this paper rely on
the following lemma and its corollary. 
The lemma is well-known, but we could not find a reference and saw no
harm in recording the proof.   We thank Jim Humphreys for simplifying our earlier proof.

\pagebreak
\begin{lem} \label{momma_lemma}
Let $x, y \in W$. Then the following four conditions are equivalent:~
\begin{enumerate}
\item[(i)] $N(x) \subseteq N(yx)$
\item[(ii)] $x^{-1} N(y) \subseteq \posroots$
\item[(iii)] $N(yx) = N(x) \cup x^{-1} N(y)$
\item[(iv)] $l(yx) = l(y) + l(x)$
\end{enumerate}
\end{lem}

\begin{proof}
Certainly (iii) implies (ii) by definition.
Then
\begin{equation} \label{equivalence}
N(x) \subseteq N(yx) \text{  if and only if  }  
x^{-1} N(y) \subseteq \posroots
\end{equation}
since $\al \in x^{-1}N(y)$ with $\al \prec 0$ if and only if 
$-\al \in N(x)$ and $-\al \notin N(yx)$. 
This shows the equivalence of (i) and (ii). 

It is easy to check from the definitions that 
\begin{equation}\label{containment}
N(yx) \subseteq  N(x) \cup  x^{-1}N(y) 
\end{equation}
where the union is a disjoint union.
Also
\begin{equation} \label{lengths}
l(yx) \leq l(x) + l(y)
\end{equation} 
which implies that 
$$|N(yx)| \leq |N(x)| + |N(y)| = |N(x)| + |x^{-1}N(y)|.$$

Hence equality in Equation (\ref{lengths}) is equivalent to equality in Equation (\ref{containment}),
which shows that (iii) and (iv) are equivalent.
Finally, it is clear that 
$$x^{-1}N(y) \cap \posroots \subseteq N(yx),$$
so if (i) holds so does (iii) by Equation (\ref{equivalence}) and Equation (\ref{containment}).
\end{proof}

The next corollary follows from the lemma by taking a reduced expression for 
$y = w x^{-1} = s_{\beta} \cdots s_{\al}$ where $\al, \beta \in \Pi$.

\begin{cor} \label{momma_cor_2}
Given $x, w \in W$.  If $N(x) \subsetneq N(w)$,
then there exists $\al, \beta \in \Pi$ 
so that 
\begin{eqnarray*}
N(x)  \subsetneq N(s_{\al}x) \subseteq N(w) & \text{ and }\\
N(x)  \subseteq  N(s_{\beta}w) \subsetneq N(w) &
\end{eqnarray*}
\end{cor}

We conclude the section with a well-known lemma 
related to the minimal length coset representatives of a parabolic subgroup of $W$.
We include a proof since we refer to the proof in what follows.


Let $\ro' \subset \ro$ be a parabolic subsystem. 
In other words, $\ro'$ has a basis of simple roots which is contained in 
the simple roots $\Pi$ of $\ro$.
Let $\ro^1 = \ro^+ - \ro'^+$ and let $W' \subset W$ be the Weyl group of $\ro'$.

\begin{lem} \label{coset_reps}
The set $C:= \{ x \in W \ | \ N(x) \subseteq \ro^1 \}$ 
is a set of distinct coset representatives for $W'$ in $W$.
\end{lem}

\begin{proof}
Take $w \in W$. 
The intersection $N(w) \cap \ro'$ is of Weyl-type for $\ro'$.  Thus
$$N(w) \cap \ro' = N(x)$$
for some $x \in W'$, where $N(x)$ is
the same whether computed in $W'$ or $W$.
Now $N(x) \subseteq N(w)$, so by Lemma \ref{momma_lemma},
we have that $y = wx^{-1}$ satisfies
$$x^{-1}N(y) = N(w) - N(x).$$
Hence $x^{-1}N(y) \subseteq \ro^1$ by the definition of $x$.
Consequently $N(y) \subseteq \ro^1$ as $x \in W'$ and $W'$
preserves $\ro^1$.
Certainly, $wW' = yxW' = yW'$, which shows that the elements of $C$ are a set 
of coset representatives.

They must be distinct representatives.
Indeed, suppose $y = z x$ for $y,z \in C$ and some $x \in W'$.
Then $x^{-1}N(z) \subseteq \ro^1$ since $W'$ preserves $\ro^1$.
Then Lemma~\ref{momma_lemma}~ implies that 
$N(x) \subseteq N(y) \subseteq \ro^1$, forcing $x =1$ since also
$N(x) \subseteq \ro'$ and so $y=z$.


\end{proof}

\section{Some new results for general root systems}

The first new result of this paper is a generalization
of the fact that every subset of $\posroots$ of Weyl-type is of the form $N(w)$ for some $w \in W$
\cite{kostant:borel}.
On the one hand, it is easy to see that that
$N(w) \cap \ideal^c$  for $w \in W$ is always of Weyl-type 
in $\ideal^c$.  But the converse is also true.  
Namely, every subset $S$ of $\ideal^c$ of Weyl-type
is of the form $N(w) \cap \ideal^c$ for some $w \in W$.
That is, 

\begin{prop}
Let $S \in {\mathcal W}^{\ideal}$.  There exists $w \in W$ such that
$S = N(w) \cap \ideal^c$.
\end{prop}

\begin{proof}
This is equivalent to showing that there exists $T \in {\mathcal W}^{\emptyset}$
for which $S = T \cap \ideal^c$ since we already know the result is true when $\ideal = \emptyset$.
In fact, it is enough to prove that given an ideal $\ideal_1$ where 
$\ideal_1^c = \ideal^c \cup \{ \delta \}$ for $\delta \in \posroots$ we can find 
$T \in {\mathcal W}^{\ideal_1}$ where $S = T \cap \ideal^c$.  Then the result
would follow by induction as there is always a sequence
$$\ideal^c = \ideal_0^c \subset  \ideal_1^c \subset \dots \subset  \ideal_r^c = \posroots,$$
such that 
$\ideal_i$ is an ideal and 
$|\ideal_{i+1}^c| = |\ideal_{i}^c| + 1$. 

There are four possible situations given the above setup:
\begin{enumerate}
\item $S, \ S \cup \{ \delta \} \in {\mathcal W}^{\ideal_1}$
\item $S \in {\mathcal W}^{\ideal_1}$, \ $S \cup \{ \delta \}  \notin {\mathcal W}^{\ideal_1}$
\item $S \notin {\mathcal W}^{\ideal_1}$, \ $S \cup \{ \delta \}  \in {\mathcal W}^{\ideal_1}$
\item $S, \  S \cup \{ \delta \} \notin {\mathcal W}^{\ideal_1}$
\end{enumerate}
Assuming one of the first three possibilities arises, one of 
$S$ or $S \cup \{ \delta \}$ (or in the first case both) would suffice for $T$.
The proposition is therefore equivalent to the last possibility never occurring.

The last possibility would only occur if there exists $\al, \beta \in S$
and $\al', \beta' \in \ideal^c - S$ for which 
$\al + \beta = \delta = \al' + \beta'$.
If so,
then 
$$ -\al' + \al + \beta = \beta'.$$
By Lemma 3.2 in \cite{sommers:b_stable}, either 
$-\al' + \al$ or $-\al' + \beta$ lies in $\ro \cup \{ 0 \}$.
Since $\al, \beta \in S$ and $\al' \notin S$, neither 
$-\al' + \al$ nor $-\al' + \beta$ can be zero.
Without loss of generality we take
$-\al' + \al \in \ro$, and by possibly interchanging the roles of $\al$ and $\al'$ in what
follows, we may assume $-\al' + \al \in \posroots$.

On the one hand, $(-\al' + \al) + \al' = \al$.
Now since $\ideal$ is an ideal and $\al \in \ideal^c$
and $-\al' + \al \prec \al$, then $-\al' + \al \in \ideal^c$.
It follows that $-\al' + \al \in S$, since $\al' \notin S$ and $\al \in S$
and otherwise $S^c = \ideal^c - S$ would not be $\ideal^c$-closed.

On the other hand, $-\al' + \al = \beta' - \beta \in \ideal^c$.
Clearly, $(\beta' - \beta) + \beta = \beta'$.  Then 
$\beta' - \beta \notin S$, since $\beta \in S$ and $\beta' \notin S$
and otherwise $S$ would not be $\ideal^c$-closed.
This contradicts the fact that $-\al' + \al = \beta' - \beta \in S$ from 
the previous paragraph.

It follows that $T$ can be chosen to be 
one of $S$ or $S \cup \{ \delta \}$ (or possibly both), proving the proposition.
\end{proof}

Although there is not always a unique $w \in W$ satisfying the hypotheses of the
proposition, there is a unique $w$ with the property that
$N(w)$ is contained in $N(w')$ for any other $w' \in W$ 
satisfying the hypotheses of the
proposition.
More generally, 

\begin{prop} \label{minimal}
Let $\ideal' \subseteq \ideal$ be ideals.
Given $S \in {\mathcal W}^{\ideal}$, there exists 
$T \in {\mathcal W}^{\ideal'}$ with the property
that $S = T  \cap \ideal^c$ and if  
${\widehat T} \in {\mathcal W}^{\ideal'}$ satisfies
$S \subseteq {\widehat T}$, then $T \subseteq {\widehat T}$.  Such a $T$ is clearly unique.
\end{prop}

\begin{proof}
The proof is by induction on the difference in cardinalities
$l = |\ideal| - |\ideal'|$.  The case $l=0$ is trivial, with $T=S$.
If $l>0$, pick an ideal $\ideal_1$ 
such that 
$$\ideal' \subsetneq \ideal_1 \subseteq \ideal$$
where $\ideal'^c = \ideal_1^c \cup \{ \delta \}$
for some $\delta \in \posroots$.
By induction there exists 
$T_1 \in {\mathcal W}^{\ideal_1}$ satisfying the hypotheses
of the proposition with respect to $\ideal_1 \subseteq \ideal$ 
and $S \in {\mathcal W}^{\ideal}$.  

Set $T= T_1$ if $T_1 \in {\mathcal W}^{\ideal'}$, or else
set $T= T_1 \cup \{ \delta \}$ if $T_1 \notin {\mathcal W}^{\ideal'}$.
In either case the proof of the previous proposition ensures that
$T \in {\mathcal W}^{\ideal'}$.

Now suppose that ${\widehat T} \in {\mathcal W}^{\ideal'}$ satisfies
$S \subseteq {\widehat T}$.
Since $\ideal_1^c  \subset \ideal'^c$, it is clear that 
${\widehat T} \cap \ideal_1^c \in {\mathcal W}^{\ideal_1}$.
Then the minimal property for $T_1$ gives that 
$$T_1  \subseteq {\widehat T} \cap \ideal_1^c.$$ 
Hence $T_1 \subseteq {\widehat T}$.  
We deduce that 
$T \subseteq {\widehat T}$. 
Indeed, $T =  T_1 \cup \{ \delta \}$ only when 
there exist $\al, \beta \in T_1$ with $\al + \beta = \delta$.
Since $T_1 \subseteq {\widehat T}$ and ${\widehat T}$ is $\ideal'^c$-closed,
we must have $\delta \in {\widehat T}$.
\end{proof}

There is a nice characterization of the 
$w \in W$ for which $T = N(w)$ satisfies
the minimal condition of Proposition \ref{minimal} when $\ideal' = \emptyset$.

\begin{prop} \label{condition}
Given $S \in {\mathcal W}^{\ideal}$
there is a unique $w \in W$ satisfying
both $S = N(w) \cap \ideal^c$ and 
\begin{equation} \label{nonempty}
w^{-1}(-\Pi) \cap \posroots \subseteq \ideal^c.
\end{equation}
Furthermore, $N(w) \subseteq N(x)$ for any $x \in W$ with
$S \subseteq N(x)$ as in Proposition \ref{minimal}.
\end{prop}

\begin{proof}
Let $w$ be such that $T = N(w)$ 
satisfies the minimal property from Proposition \ref{minimal}
for $\ideal'= \emptyset$ and $S$.

Suppose there exists a simple root $\al \in \Pi$ for which
$$w^{-1}(-\al) \in \posroots - \ideal^c.$$
Consider $x = s_{\al} w$.  Then 
$$w^{-1}(-\al) = x^{-1}(\al) \in \posroots$$ 
and so Lemma \ref{momma_lemma} implies that 
$$N(w) = N(x) \cup \{ w^{-1}(-\al) \}.$$
But $w^{-1}(-\al) \notin \ideal^c$ and
therefore 
$$N(x) \cap \ideal^c = N(w) \cap \ideal^c,$$
contradicting the minimal property of $N(w)$.  Hence 
we must have
$$w^{-1}(-\Pi) \cap \posroots \subseteq \ideal^c.$$

For the uniqueness, take $y \in W$ with
$S = N(y) \cap \ideal^c$ and $y \neq w$. 
Then $N(w) \subsetneq N(y)$ by Proposition \ref{minimal}.
By Corollary \ref{momma_cor_2}
there exists $\al \in \Pi$ such that 
$$N(w) \subseteq N(s_{\al} y) \subsetneq N(y)$$
and this implies that $N(y)= N(s_{\al} y) \cup \{ y^{-1}(-\al) \}$.
It follows that $y^{-1}(-\al) \notin \ideal^c$ since
$$N(y) \cap \ideal^c = N(w) \cap \ideal^c$$ and hence
$y^{-1}(-\Pi) \cap \posroots \not\subset \ideal^c$.
\end{proof}

\section{Further results for types $A_n, B_n, C_n$}

In this section we explore some properties which are particular to 
types $A_n, B_n, C_n$ and which are used in the proof of Theorem \ref{thm1}.

We label the simple roots in types $B_n$ and $C_n$ so that $\al_n$ is the only simple root of its length.  
For each root system of type $X_n$, we embed
$X_{n-1}$ in $X_{n}$ via the simple roots $\al_2, \dots, \al_n$.
Denote by $\ro_{n-1}$ the roots of $X_{n-1}$ and 
let 
$$\ro^1 = \posroots - \posroots_{n-1}.$$
Let $W' \subset W$ be the Weyl group of $X_{n-1}$. 
One of the key facts about these root systems
is that $\ro^1$ is linearly ordered under $\prec$ and (thus) there
is only root of each height from $1$ to the largest height ($n$ 
in type $A_n$ and $2n-1$ in types $B_n, C_n$).

More precisely,

\begin{lem} \label{abc_lemma}
Given $\al, \beta \in \ro^1$ with $\al \prec \beta$, 
we always have
\begin{equation} \label{key_one}
\beta - b  \al = c  \gamma
\end{equation}
for some $\gamma \in \ro_{n-1}$ and 
some $b, c \in \{ 1, 2 \}$.
\end{lem}

\begin{proof}
The roots in $\ro^1$ (in a standard basis)
are
\begin{eqnarray*}
&\{ e_1 - e_j  \  | \ j = 2, 3, \dots, n \} & \text{ for } A_{n-1}  \\
&\{ e_1 \pm e_j  \  | \ j = 2, 3, \dots, n\} \cup \{e_1\}  & \text{ for } B_n \\
&\{ e_1 \pm e_j  \  | \ j = 2, 3, \dots, n\} \cup \{2e_1\}  & \text{ for } C_n  \\
\end{eqnarray*}
These roots are ordered as follows:
$$e_1 - e_2 \prec \dots \prec e_1 - e_n \prec e_1 
\prec e_1+e_n \prec \dots \prec e_1 +e_2 \prec 2e_1,$$
whenever the given root is present in the appropriate root system.

It is easy to see in type $A_{n-1}$ that $\beta - \al \in \ro^+_{n-2}$ (as desired, 
given the shift in subscript).

In type $B_n$, $\beta - \al \in  \ro^+_{n-1}$, unless $\beta = e_1 + e_j$ 
and $\al = e_1 - e_j$, in which case $\beta - \al = 2 \gamma$
for some $\gamma \in \ro^+_{n-1}$.

In type $C_n$, we have $\beta - \al \in \ro^+_{n-1}$, except when $\beta = 2e_1$,
the highest root.  In that case, we can say that $\beta - 2 \al \in \ro_{n-1}$. 
\end{proof}

\begin{lem} \label{key_lemma}
In types $A_n, B_n, C_n$, 
for each $k \in \{0, 1, \dots, |\ro^1| \}$,
there is a unique element $x \in W$ satisfying 
$N(x) \subseteq \ro^1$ and $|N(x)| = k$. 
\end{lem}


\begin{proof}
Certainly there exists a unique $w^1$
with $N(w^1) = \ro^1$ since $\ro^1 \in {\mathcal W}^{\emptyset}$.
Explicitly, $w^1$ is the product of the long element of $W$ and the long element of $W'$,
using, for example, Lemma \ref{momma_lemma}.
It follows that there exists at least one $x \in W$ satisfying 
$N(x) \subseteq \ro^1$ and $|N(x)| = k$
from Corollary \ref{momma_cor_2} by taking any reduced expression for $w^1$.

We still need to show uniqueness.  
Let $x \in W$ satisfy $N(x) \subsetneq \ro^1$ and 
assume uniqueness is true when $k > |N(x)|$.
Since $N(x) \subsetneq N(w^1)$, 
Corollary \ref{momma_cor_2} implies the existence of $\al \in \Pi$
such that 
$$N(x) \subsetneq N(s_{\al}x) \subseteq N(w^1),$$
where $N(s_{\al}x) = N(x) \cup \{ x^{-1}(\al) \}$.
Since $N(s_{\al}x) \subseteq N(w^1) = \ro^1$, we have
$$x^{-1}(\al) \in \ro^1.$$ 

We claim that $\al$ is unique.  Indeed, assume that 
$\beta \in \Pi$ and $x^{-1}(\beta) \in \ro^1$ with $\beta \neq \al$.
Without loss of generality,
$x^{-1}(\al) \prec x^{-1}(\beta)$.  By Lemma \ref{abc_lemma},
$x^{-1}(\beta) - b x^{-1}(\al) = c \gamma$ for some $b,c$ and $\gamma \in \ro_{n-1}$.
Applying $x$ to both sides yields
$\beta - b  \al = c  x(\gamma)$.  This is impossible
since the right side is a linear combination of simple roots with
either all positive or all negative coefficients, whereas the left side
is a combination of two simple roots whose coefficients have opposite signs.
Since $s_{\al}x$ is unique by induction, $x$ is unique and the result follows.

\end{proof}

\section{Proof of Theorem \ref{thm1} for types $A_n, B_n, C_n$} \label{main_proof}


Assume the factorization is true for $\ro_{n-1}$.
Clearly $\ideal^c \cap  \ro_{n-1}$ is equal to $\ideal'^c$ for 
some ideal $\ideal'$ for $\ro_{n-1}$.
Let $m_1, \dots, m_{k-1}$ be the $\ideal'$-exponents.  

As noted in the previous section, the roots of 
$\ro^1 = \posroots - \posroots_{n-1}$ 
are linearly ordered and so contain one root of each height.
It follows that
$\ideal^c \cap \ro^1$ contains one root of height
$1, 2, \dots, m_k$ for some natural number $m_k$, 
and that 
the $\ideal$-exponents are $m_k$ together with the $\ideal'$-exponents
$m_1, m_2, \dots, m_{k-1}$ (in this indexing, $m_k$ need not be 
the largest exponent).

Now it is certainly true that if $S \in {\mathcal W}^{\ideal}$, then 
$S \cap \ro_{n-1} \in {\mathcal W}^{\ideal'}$.
We have a strong converse which holds in types $A_n ,B_n, C_n$:

\begin{lem}
Let $S' \in {\mathcal W}^{\ideal'}$. 
For each $j \in \{ 0, 1, \dots, |\ro^1| \}$, there exists
a unique $S \in {\mathcal W}^{\ideal}$
with $S \cap \ro_{n-1} = S'$ and $|S \cap \ro^1| = j$.
\end{lem}

\begin{proof}
Let $w' \in W'$  be the unique element with 
the property that $N(w') \cap \ideal'^c = S'$
and $N(w') \subseteq N(x)$ for any $x \in W$ with $N(x) \cap \ideal'^c = S'$ 
as in Proposition \ref{minimal}.  

Let $\{ x_0=1, x_1, \dots, x_i, \dots \}$ be the elements from Lemma \ref{key_lemma}
with $|N(x_i)|=i$.  Note that 
\begin{equation*} 
N(x_i) \subsetneq N(x_{i+1})
\end{equation*}
from the existence part of the 
proof of that lemma.

Consider the elements $x_i w'$.
As in the proof of Lemma \ref{coset_reps},
\begin{equation} \label{yeah} 
N(x_i w')  = N(w') \cup w'^{-1} N(x_i)
\end{equation}
where $w'^{-1} N(x_i) \subseteq \ro^1$.
It follows that $N(x_{i} w') \subsetneq N(x_{i+1} w')$
and the two sets differs by a single element of $\ro^1$.

Next, consider the intersection 
$$N(x_i w') \cap \ideal^c \cap \ro^1.$$
This intersection is empty for $i=0$ and has $m_k$ elements
when $i=n$ in type $A_n$ and when $i=2n-1$ in types $B_n$ and $C_n$.
From the previous paragraph, 
we know that 
$N(x_{i+1} w') \cap \ideal^c \cap \ro^1$ and 
$N(x_i w') \cap \ideal^c \cap \ro^1$ can differ by at most one
element.
Consequently, for some $i$ we have that $S:=N(x_i w') \cap \ideal^c$
satisfies $S \cap \ro_{n-1} = S'$ and  $|S \cap \ro^1|=j$
has the desired cardinality.
This gives existence.  It would also give uniqueness if we
knew that every $S$ is of the form $N(x_i w') \cap \ideal^c$
for some $i$.

To that end, suppose that 
$S \in {\mathcal W}^{\ideal}$ and $S \cap \ro_{n-1} = S'$.
Let $w \in W$ be the unique element with 
the property that $N(w) \cap \ideal^c = S$
and $N(w) \subseteq N(x)$ for any $x \in W$ with $S \subseteq N(x)$
as in Proposition \ref{minimal}.

Write $w = x_i y$ for $y \in W'$ by Lemma \ref{coset_reps}.
By Equation (\ref{yeah}) and the line following it, 
$N(y) = N(w) \cap \ro_{n-1}$. 
The latter contains $S'$ by the definition of $w$.
Thus there exists $i'$ so that 
$$S \subseteq N(x_{i'} w'),$$
simply by taking the largest possible value of $i'$.
By Proposition \ref{minimal}, we get $N(x_i y) \subseteq N(x_{i'} w')$
and thus $N(y) \subseteq N(w')$ after intersecting with $\ro_{n-1}$. 
Now Proposition \ref{minimal} applied to $w'$ gives 
the equality $N(y)= N(w')$ and so $w' = y$.
It follows that $w = x_i w'$, showing uniqueness of $S$.
\end{proof}

\noindent {\bf Proof of Theorem \ref{thm1}: }
By the previous lemma, if we consider the sum 
$\sum_{S} t^{|S|}$ over all $S \in {\mathcal W}^{\ideal}$ 
with $S \cap \ro_{n-1} = S'$ for some $S' \in {\mathcal W}^{\ideal'}$, then the sum
equals
$$t^{|S'|}(1 + t +t^2 + \dots + t^{m_k}).$$ 

Thus
\begin{eqnarray*}
\sum_{S \in {\mathcal W}^{\ideal}} t^{|S|} = 
\sum_{S' \in {\mathcal W}^{\ideal'}} t^{|S'|}(1 + t +t^2 + \dots + t^{m_k}) = \\
(1 + t +t^2 + \dots + t^{m_k}) \prod_{i=1}^{k-1} (1 +t + t^2 + \dots + t^{m_i}),
\end{eqnarray*}
where the last step is by induction.  This completes the proof
of the theorem in types $A_n, B_n, C_n$.  This proof also works in type $G_2$.  In types $F_4$ and $E_6$, 
the theorem was checked on a computer, running through all possible ideals.  There are $105$ ideals in $F_4$ and
$833$ of them in $E_6$.  
 
\section{A uniform proof for the penultimate ideal} \label{uniform}

The goal for this section is to prove Theorem \ref{thm1} uniformly 
when $\ideal = \{ \theta \}$, where $\theta$ is the highest root of $\posroots$.

The next result is a special case of Theorem 2.8 in \cite{macdonald:poincare}.
For $\gamma$ in $X$ (the weight lattice), 
let $e^{\gamma}$ denote the corresponding element of the
group algebra $\Z[X]$ of $X$.
\begin{prop}  \cite{macdonald:poincare}
Let $R \subset \posroots$ be any subset of the positive roots.
The following identity holds: 
\begin{equation} \label{macky}
\sum_{w \in W} \prod_{\al \in R} \frac{1 - t e^{w \al}}{1 - e^{w \al}} 
=
\sum_{w \in W} t^{|N(w) \cap R|}  
\end{equation}
\end{prop}

\begin{proof}
In Theorem 2.8 of 
\cite{macdonald:poincare}, 
set $u_{\al} = t$ if $\al \in R$ and set $u_{\al} = 1$ if $\al \not\in R$. 
\end{proof}

When $R = \ideal^c$ for some ideal $\ideal$, the right-side of 
(\ref{macky}) is the Poincar{\'e} polynomial (after replacing $t$ by $t^{2}$) 
of a regular semisimple Hessenberg variety
(see the next section) by work of 
\cite{MPS}.  In that case, the identity can be proven by a fixed-point formula
as in \cite{macdonald:poincare} since these Hessenberg 
varieties are smooth and projective (generalizing
the role of the flag variety in Macdonald's fixed-point formula proof).

We now use this 
identity to prove Theorem \ref{thm1} uniformly for any root system 
when $\ideal = \{ \theta \}$.



\begin{thm} \label{uniform1}
In all types when $\ideal = \{ \theta \}$, 
$$\sum_{S \in {\mathcal W}^{\ideal}} t^{|S|} 
= \prod_{i=1}^{n} (1 +t + t^2 + \dots + t^{m^{\ideal}_i})$$
where the $m^{\ideal}_i$ are the exponents of $\ideal$.
\end{thm}
\begin{proof}

In Equation (\ref{macky}), put $R = \ideal^c = \posroots - \{ \theta \}$ 
and specialize $e^{\al}$ to 
$t^{\hgt(\al)}$ as in \cite{macdonald:poincare}.  Then we have
$$\sum_{w \in W} \prod_{\al \in \ideal^c} \frac{1 - t^{ \hgt(w \al) + 1 }}{1 - t^{ \hgt(w \al)}} 
=
\sum_{w \in W} t^{|N(w) \cap \ideal^c|}  .$$


We will break apart the sum on the left side into two parts, according to whether
$w \in W$ satisfies Equation (\ref{nonempty}).
Let $W_{min}$ denote those elements of $W$ satisfying Equation (\ref{nonempty})
for $\ideal$.

If $w \in W_{min}$ and $w \neq 1$, then 
$w^{-1}(\al) \in \posroots$ for some $\al \in -\Pi$ since $w \neq 1$.
Then $\beta := w^{-1}(\al) \in \ideal^c$ by Equation (\ref{nonempty}).
Since $\hgt(w \beta) = -1$ , the term for $w$ vanishes in the sum.  
This leaves only the identity term as the contribution from the elements in $W_{min}$.

Therefore the sum on the left side reduces to 
\begin{equation} \label{another_eq}
\prod_{\al \in \ideal^c} \frac{1 - t^{\hgt(\al)+1}}{1 - t^{\hgt(\al)}} 
+ \sum_{w \not\in W_{min}} \prod_{\al \in \ideal^c} \frac{1 - t^{\hgt(w \al)+1}}{1 - t^{\hgt(w \al)}}
\end{equation}

The isolated product is exactly the right side of Theorem \ref{uniform1}.
To finish the proof we must, by Proposition \ref{condition}, show that the sum 
in Equation (\ref{another_eq}) is equal to 
$$\sum_{w \not\in W_{min}} t^{|N(w) \cap \ideal^c|} \ .$$ 
We divide $\ideal^c$ into two parts:
let $$\ideal^c_i = \{ \gamma \in \ideal^c \ | \ \langle \gamma, \theta^{\vee} \rangle = i \}$$
for $i=0,1$ (the only two possibilities since $\theta$ is the highest root).
The roots of $\ideal^c_0$ are the positive roots of a parabolic subsystem of $\ro$, with
corresponding Weyl group $W_{\theta}$.   
These are exactly the elements of $W$ which fix $\theta$.

Given $\gamma \in \ideal^c_1$, we have 
$s_{\theta}(-\gamma) = \theta - \gamma$, which is a positive root. 
Hence also $s_{\theta}(-\gamma) \in \ideal^c_1$
and $\gamma + s_{\theta}(-\gamma) = \theta$.
This shows that elements of $\ideal^c_1$ come in pairs which sum up to $\theta$.

Take $w \not\in W_{min}$.
Then $w \theta \in -\Pi$ and thus $\hgt(w \theta)=-1$.  
Suppose that $\al + \beta = \theta$ for $\al, \beta \in \ideal^c_1$.
Then
\begin{equation} \label{pairs}
\hgt(w \al) + \hgt(w \beta) = -1.
\end{equation}
Consequently exactly one of $\al$ and $\beta$ belongs
to $N(w) \cap \ideal^c$
and thus
$$|N(w) \cap \ideal^c_1| = \frac{1}{2} |\ideal^c_1|$$
for all $w \not\in W_{min}$.

On the other hand, 
the identity 
$$\frac{1 - t^{a+1}}{1-t^a} \cdot \frac{1 - t^{-a}}{1-t^{-a-1}}
= t$$
and Equation (\ref{pairs}),
imply that
$$\prod_{\al \in \ideal^c_1} \frac{1 - t^{\hgt(w \al)+1}}{1 - t^{\hgt(w \al)}}
= t^{\frac{1}{2} |\ideal^c_1|}$$
for $w \not\in W_{\min}$.

Therefore the proof will be completed
if we can show that 
$$\sum_{w \not\in W_{min}} \prod_{\al \in \ideal_0^c} \frac{1 - t^{\hgt(w \al)+1}}{1 - t^{\hgt(w \al)}} = \sum_{w \not\in W_{min}} t^{|N(w) \cap \ideal^c_0|} \ .$$

We can do this by using Equation (\ref{macky}) for the Weyl group $W_{\theta}$.
First we observe that the action of $W_{\theta}$ on $W$ preserves $W_{min}$
and so $W_{min}$ and its complement are a union of left cosets of $W_{\theta}$.

Pick $x \not\in  W_{min}$.
Then 
\begin{equation} \label{coset}
\sum_{w \in x W_{\theta}}  \prod_{\al \in \ideal_0^c} \frac{1 - t^{\hgt(w \al)+1}}{1 - t^{\hgt(w \al)}}
=
\sum_{y \in W_{\theta}} \prod_{\al \in \ideal^c_0} \frac{1 - t^{\hgt(x y (\al))+1}}{1 - t^{\hgt(x y (\al))}}
\end{equation}

In Equation (\ref{macky}) applied now for the case of the Weyl group $W_{\theta}$,
set $e^{\al} = t^{\hgt(x \al)}$ where the height 
is still calculated with respect to $W$ (since $x \al$ need
not belong to the root system of $W_{\theta}$). 
The positive roots for $W_{\theta}$ are $\ideal_0^c$, so 
Equation (\ref{coset}) is equal to 
$$\sum_{y \in W_{\theta}} t^{|N(y) \cap \ideal_0^c|}$$
which is the same thing as
$$\sum_{w \in x W_{\theta}} t^{|N(w) \cap \ideal_0^c|}$$  
for any $x \in W$.
Indeed, we may choose $x$ such that $N(x) \subseteq \posroots - \ideal_0^c$
by Lemma \ref{coset_reps}.  Then Lemma \ref{momma_lemma} 
implies that for $w_1, w_2 \in W_{\theta}$ that
$N(xw_1) \cap \ideal_0^c = N(xw_2) \cap \ideal_0^c$ if and only
if $w_1 = w_2$, which is what we needed.  
The proof is completed since $W_{min}$ is a union of left cosets of $W_{\theta}$.

\end{proof}

\section{Poincar{\' e} polynomials of regular nilpotent Hessenberg varieties}

The combinatorial Poincar{\' e} polynomials from Theorem \ref{thm1} 
should arise as the actual topological Poincar{\' e} polynomials 
of certain projective subvarieties of the flag variety.
This section defines these subvarieties, called 
regular nilpotent Hessenberg varieties, and lists some of their main
properties.

Write $B^-$ for the Borel subgroup opposite to $B$ and $\mathfrak{b}^-$
for its Lie algebra.


Given an ideal $\ideal$, we define 
\[H_{\ideal} = \mathfrak{b}^- \oplus \bigoplus_{\alpha \in \ideal^c} \mathfrak{g}_{\alpha},\]
where $\mathfrak{g}_{\alpha} \subset \g$ is the $\al$-weight space. 
The subspace $H_{\ideal}$ is stable for the action of $B^-$ and 
it is easy to see that each subspace with this property is of the above
form.  Such a subspace is called a Hessenberg space. 



Fix an
element $X \in \mathfrak{g}$ and a Hessenberg space $H = H_{\ideal}$.  
The Hessenberg variety $\H_{X,H}$
is the subvariety of the flag variety $\H = G/B^{-}$ defined by
\[\H_{X,H} = \{gB^{-} \ | \ \Ad(g^{-1})(X) \in H\}.\]
This is a closed subvariety of $\H$ and hence is projective.  In general
a Hessenberg variety is not smooth.
Hessenberg varieties were first defined in \cite{MPS}.

When $\ideal = \posroots$ and thus $H = \bo^{-}$, 
the Hessenberg variety reduces to a Springer variety, a well-studied and important
object in representation theory.
At the other end of the spectrum, when $\ideal = 0$
then the Hessenberg variety is the whole flag variety, independent of $X$.
In between, when $\ideal^c$ is the set of simple roots
and $X$ is regular nilpotent, 
the Hessenberg variety is called the Peterson variety and 
has been used to give geometric constructions for the
quantum cohomology of the flag variety (see \cite{Ko2}, \cite{R}).
Other Hessenberg varieties have been used in \cite{GKM} to give a partial
proof of the fundamental lemma of the Langlands program.

The following proposition about $\H_{X,H}$ follows from work in \cite{T}.
Let $B^{-} w B^{-} \subset \H$ be the Schubert cell containing
the point $wB^{-}$ where $w \in W$.  
Here we do not distinguish between $w \in W$ and a representative of $w$
in $G$.

\begin{prop} \label{tymoczko thesis}
Let $G$ be of classical type and let $X$ be a sum of negative simple
root vectors.
Let $C_w := B^{-} w B^{-} \cap \H_{X,H_{\ideal}}$.  Then $C_w$ is 
non-empty if and only if $w$ satisfies Equation (\ref{nonempty}) for $\ideal$.
If $C_w$ is non-empty, then it is an affine space of dimension
$|N(w) \cap \ideal^c|$.

This yields an affine paving of $\H_{X,H_{\ideal}}$.
\end{prop}

This proposition allows us to demonstrate that the Poincar{\' e} polynomials 
for regular nilpotent Hessenberg varieties are the 
polynomials which arise in Theorem \ref{thm1}.

\begin{thm}
Let $\mbox{P}_{\ideal}(t)$ denote the Poincar{\' e} polynomial of the 
Hessenberg variety $\H_{X,H_{\ideal}}$ for $X$ a regular nilpotent element.  
In types $A_n$, $B_n$, and $C_n$ this can be factored
\[\mbox{P}_{\ideal}(t^{\frac{1}{2}}) 
= \prod_{i=1}^k (1 + t + \ldots + t^{m_i^{\I}}).\]
\end{thm}

\begin{proof}
Proposition \ref{condition} shows that for each 
$S \in \mathcal{W}^{\I}$ there is exactly one 
$w \in W$ satisfying Equation (\ref{nonempty}) 
and $S = N(w) \cap \I^c$.
Proposition \ref{tymoczko thesis} shows that the
dimension of the affine cell $C_w$ is $|N(w) \cap \I^c|=|S|$
and that $C_w$ is empty if $w$ does not satisfy Equation (\ref{nonempty}).   
The proof follows from the fact
that these cells give an affine paving of the variety together with Theorem \ref{thm1}.
\end{proof}





We should mention that if one knew that Proposition \ref{tymoczko thesis}
were true in all types 
then a result of Peterson announced in \cite{BC} would be equivalent to 
Theorem \ref{thm1}.
Unfortunately Peterson's proof is not given.  One could imagine
something along the lines of \cite{carrell:kostant-macdonald}
if the regular nilpotent Hessenberg varieties were smooth.  But this is not
the case, already in type $A_2$.

\section{Hyperplane arrangement defined by an ideal}

The second venue where the $\ideal$-exponents arise is in the context of hyperplane arrangements.
Let $V := \colat \otimes \mathbf{R}$ be the ambient vector space containing the coroot lattice
$\colat$.
For each $\al \in \posroots$, let $H_{\al} \subset V$ be the hyperplane 
$$H_{\al} = \{ v \in V \ | \ \langle \al, v \rangle = 0 \}.$$
We are interested in the hyperplane arrangement in $V$
given by the hyperplanes $H_{\al}$ where $\al \in \ideal^c$.
We will denote this arrangement by $\hyp_{\ideal}$ and call it an arrangement of ideal-type
in $\ro_n$.  

In general, 
given a hyperplane arrangement,
one is interested in whether the arrangement is free and if so,
what are the roots of its characteristic polynomial, which are also called
exponents \cite{orlik-terao}. 

We briefly recall the basic definitions and theorems about hyperplane arrangements
from Chapters 2 and 4 of \cite{orlik-terao}.
Let $\A$ be an arrangement of hyperplanes in the $\mathbb R$-vector space $V$.
Let $S(V^*)$ denote the symmetric algebra on the dual space $V^*$ of $V$.
Given $H \in \A$, let $\al_H \in V^*$ be a linear functional vanishing on $H$.
Set 
$$Q(\A) = \prod_{H \in \A} \al_H.$$
Let $D(\A)$ denote the $\mathbb R$-linear 
derivations of $S(V^*)$ which preserve the ideal generated
by $Q(\A)$.
Then the hyperplane arrangement $\A$ is said to be free if $D(\A)$ is a free $S(V^*)$-module. 

Let $L=L(\A)$ denote the set of nonempty intersections of the elements 
of $\A$.  This is a poset with a partial order given by reverse inclusion, with minimal element
$V$.  Define a function $\mu$ on $L$ as follows.  Set $\mu(V)=1$ and define
$\mu(X)$ recursively for $X \in L$ by the formula
$$\mu(X) = - \sum_{V \leq Z < X} \mu(Z),$$
where the sum is over all $Z \in L$ with $V \leq Z < X$.
Then the characteristic polynomial $\chi(\A, t)$ of $\A$ is defined
as
$$\chi(\A, t) = \sum_{X \in L} \mu(X)t^{\text{dim}(X)}.$$

The factorization result of Terao (see Theorem 4.137 in \cite{orlik-terao}) states
that if 
$\A$ is free, then all of the roots of $\chi(\A, t)$ are nonnegative integers, called 
the exponents of $\A$.  These exponents coincide with the polynomial degrees of a set of 
homogeneous generators of $D(\A)$ as an $S(V^*)$-module. 

There is another key property of hyperplane arrangements. 
Given $H_0 \in \A$, 
let $\A'$ denote the arrangement in $V$ obtained 
by omitting the hyperplane $H_0$ from $\A$.
This is the deleted arrangement given by $H_0$.
Let $\A''$ denote the arrangement in $H_0$ given by the nonempty intersections $H \cap H_0$
for $H \in \A$ with $H \neq H_0$.
This is the restricted arrangement given by $H_0$.
The three arrangements $(\A, \A', \A'')$
is called a triple of arrangements.

We will use the following direction
of the Addition-Deletion Theorem of Terao in what follows 
(see Theorem 4.51 in \cite{orlik-terao}).
If $\A'$ is free with exponents $b_1, \ldots, b_{k-1}, b_k-1$
and if $\A''$ is free with exponents $b_1, \ldots, b_{k-1}$,
then $\A$ is free with exponents $b_1, \ldots, b_k$.
If we wished only to know about the implication involving the exponents, 
this result goes back to Brylawski and Zaslavsky (see Theorem 2.56 in \cite{orlik-terao}).

\begin{thm} \label{thm2}
Except possibly in types $F_4, E_6, E_7, E_8$,
the hyperplane arrangement $\hyp_{\ideal}$ is free and its
non-zero exponents are $m^{\ideal}_1, \dots, m^{\ideal}_k$.
There are also $n-k$ exponents equal to $0$.
\end{thm}

\begin{proof}
We will show that $\hyp_{\ideal}$ is free with the desired exponents by 
using the Addition-Deletion Theorem.
First assume that $X_n$ is of type $A_n, B_n, C_n$ and assume the
result for any ideal $\ideal_1$ properly containing $\ideal$.
Furthermore assume
the result for root systems of smaller rank of these types.  The theorem is clearly
true for the base case where $\ideal = \posroots$, since the arrangement is empty. 

Let $\ideal_1$ be
the unique ideal for which 
$\ideal_1 = \ideal \cup \{ \delta \}$
where $\delta$ is the maximal root 
in $\ro^1 \cap \ideal^c$. 
Of course, if the latter intersection is empty, then 
we are already done since the arrangement is the direct product
of the one-dimensional empty arrangement and an arrangement of ideal-type
in $\ro_{n-1}$.  

Now by induction $\hyp_{\ideal_1}$ is free
and its non-zero exponents are 
$$m^{\ideal}_1, \dots, m^{\ideal}_{k}-1$$ (we do not order the exponents) 
where we have assigned $m^{\ideal}_k = \hgt(\delta)$.

Next consider the restricted arrangement defined
by $H_{\delta}$.  This is the arrangement
in $H_{\delta}$ defined by the hyperplanes 
$H_{\al} \cap H_{\delta}$ for $\al \in \ideal^c$ and $\al \neq \delta$.
We denote this restricted arrangement by $\mathcal A^{\delta}$.
For $\beta \in \ro^1$ and $\beta \prec \delta$
Equation (\ref{key_one}) says that 
$\delta - b \beta = c \gamma$ where $\gamma \in \ro_{n-1}$.
It follows that either $\gamma$ or $-\gamma$ is in $\ideal^c$ since $\delta \in \ideal^c$.  
Set $\ideal' = \ideal \cap \ro_{n-1}$.
Then the hyperplane arrangement 
$\hyp_{\ideal'}$ (defined for $\ro_{n-1}$)
is isomorphic to $\mathcal A^{\delta}$.  Indeed, 
for each $\beta \in \ro^1$ with $\beta \neq \delta$ we have
$$H_{\beta} \cap H_{\delta} = H_{\gamma} \cap H_{\delta}$$
for some $\gamma \in \ideal'^c$.
Thus the hyperplanes $H_{\gamma} \cap H_{\delta}$
for $\gamma \in \ideal'^c$ yield the distinct
hyperplanes in $\mathcal A^{\delta}$.  
 
Thus $\mathcal A^{\delta}$ is free and its non-zero exponents
are $$m^{\ideal}_1, m^{\ideal}_2, \dots, m^{\ideal}_{k-1},$$ 
since these are the exponents of $\ideal'$.

Consequently by Addition-Deletion applied to the triple
of arrangements $(\hyp_{\ideal}, \hyp_{\ideal_1}, \mathcal A^{\delta})$,
the arrangement $\hyp_{\ideal}$ 
is free and its non-zero exponents
are equal to 
$$m^{\ideal}_1, m^{\ideal}_2, \dots, m^{\ideal}_{k},$$
as desired.

Type $G_2$ is trivial.  We consider the case of type $D_n$.
Here, 
$$\ro^1 = \{ e_1 \pm e_j \ | \ 2 \leq j \leq n \},$$
using the standard notation for roots in $D_n$.
Let 
\begin{eqnarray*}
\gamma_1 = & e_1 + e_n \\
\gamma_2 = & e_1 - e_n .
\end{eqnarray*}

The above proof carries over perfectly well as long as $\gamma_1$ and $\gamma_2$
do not both belong to $\ideal^c \cap \ro^1$
since Equation (\ref{key_one}) would hold for $\al, \beta \in \ideal^c \cap \ro^1$.

Suppose that $\gamma_1, \gamma_2 \in \ideal^c$.
First, assume that both $\gamma_1$ and $\gamma_2$ are maximal elements of $\ideal^c \cap \ro^1$.
Then $\ideal^c$ consists of the $n$ non-empty sets:
\begin{equation} \begin{split} \label{roots}
& \{e_{1+i} - e_j \ | \ 2+i \leq j \leq n \} \cup 
\{e_{1+i} + e_ j \ | \ a_i \leq  j \leq n-1\} \text{ with } 1 \leq i \leq n-2  \\
& \{e_1 - e_j | \ 2 \leq j \leq n \} \\
& \{e_j +  e_n | \ 1 \leq j \leq n-1 \} 
\end{split} \end{equation}
where $a_i$ is some natural number satisfying $2+i \leq a_i \leq n-1$ .
The former sets contain one root of each height $1, \dots, m^{\ideal}_i$,
where $m^{\ideal}_i$ depends on $a_i$.  The latter two sets contain
one root of each height $1, \dots, n-1$ and so the ideal exponents
of $\ideal$ can be written as 
$$m^{\ideal}_1, m^{\ideal}_2, \dots, m^{\ideal}_{n-2}, n-1, n-1.$$

Let $\ideal_1 = \ideal \cup \{ \gamma_1 \}$.
By induction $\hyp_{\ideal_1}$ is free
and its non-zero exponents are 
$$m^{\ideal}_1, \dots, m^{\ideal}_{n-2}, n-1, n-2.$$

Consider the restricted arrangement  $\mathcal A^{\gamma_1}$ defined
by $H_{\gamma_1}$.  We want to show it is free
with non-zero exponents $m^{\ideal}_1, \dots, m^{\ideal}_{n-2}, n-1$.

In order to do this, consider the deleted and restricted arrangements 
of $\mathcal A^{\gamma_1}$ defined by $H_{\gamma_2}$. 
The deleted arrangement $\mathcal (A^{\gamma_1})'$ is isomorphic
to the arrangement defined by $\ideal' = \ideal \cap \ro_{n-1}$ in $\ro_{n-1}$
by the same proof as in the other classical cases.
Thus $\mathcal (A^{\gamma_1})'$  is free and its
non-zero exponents 
$$m^{\ideal}_1, \dots, m^{\ideal}_{n-2}, n-2$$
by inspection of the heights of roots in (\ref{roots}).

On the other hand, the restricted arrangement $\mathcal (A^{\gamma_1})''$
defined by $H_{\gamma_2}$ is more complicated.
This arrangement lives in $H_{\gamma_1} \cap H_{\gamma_2}$
which coincides with the intersection of the null spaces of $e_1$ and $e_n$.
The hyperplanes defining $\mathcal (A^{\gamma_1})''$
are given by the null spaces of 
$$\{e_{1+i} - e_j \ | \ 2+i \leq j < n \} \cup \{e_{1+i} \} \cup
\{e_{1+i} + e_ j \ | \ a_i \leq  j \leq n-1\},$$
for $i = 1, \dots, n-2$.
This arrangement is precisely an ideal arrangement in $B_{n-2}$
which is free with non-zero exponents equal to 
$$m^{\ideal}_1, \dots, m^{\ideal}_{n-2}.$$
By the Addition-Deletion Theorem, 
${\mathcal A}^{\gamma_1}$ is free
with non-zero exponents equal to
$$m^{\ideal}_1, \dots, m^{\ideal}_{n-2}, n-1,$$
and using the theorem a second time,
it follows that $\hyp_{\ideal}$
is free and its non-zero exponents
are $$m^{\ideal}_1, \dots, m^{\ideal}_{n-2}, n-1, n-1,$$
as desired.

Finally, consider the case where $$\delta = e_1 + e_{2n - 1 - k}$$
is the maximal element of $\ideal^c \cap \ro^1$ 
for  $k > n-1$.
In this case $\hyp_{\ideal_1}$ is free
with non-zero exponents 
$$m^{\ideal}_1, \dots, m^{\ideal}_{n-2},  n-1, k-1,$$
by induction.

It suffices to complete the proof by showing that 
the restricted arrangement
of $\hyp_{\ideal}$ defined by $H_{\delta}$ 
is isomorphic to ${\mathcal A}^{\gamma_1}$ above.
Consider the element $w$ of the Weyl group of $D_n$ given
by exchanging $e_n$ and $e_{2n - 1 - k}$ and fixing all other $e_i$.
It is not hard to check that the hyperplanes defining ${\mathcal A}^{\gamma_1}$
in $H_{\gamma_1}$ are mapped to the hyperplanes defining this restricted arrangement
in $H_{\delta}$, yielding the isomorphism and completing the proof in type $D_n$.

\end{proof}

We conclude this section by noting that a uniform proof in all types 
for the case $\ideal = \{ \theta \}$ is easy.  
On the one hand, the full Coxeter arrangment is free 
with exponents the usual exponents 
$m_1 \leq \dots \leq m_n$.
On the other hand,
the restricted arrangement for $H_{\theta}$ is free with exponents
$m_1, \dots, m_{n-1}$ by \cite{terao}.
Thus $\hyp_{\ideal}$ is free with the desired exponents.

\section{Speculation}

We speculate that the two conjectures are equivalent by a general principle.  
Namely, suppose a hyperplane arrangement $\hyp$ is free 
with exponents 
$m_1, \dots, m_n$.   
Suppose further that the arrangement is central, meaning each hyperplane contains 
the origin (what we will speculate is false otherwise).

Let $\mathcal C_{\hyp}$ denote the set of components 
of the complement $V - \cup_{H \in \hyp} H$.   
Fix one component $A \in \mathcal C_{\hyp}$.  Then for each component 
$B \in \mathcal C_{\hyp}$ we can define $l(B)$ to be the least number 
of hyperplanes needed to be crossed to move from $B$ to $A$.  

Define a polynomial
$$P_{\hyp} (t) =  \sum_{B \in \mathcal C_{\hyp}} t^{l(B)}.$$
This is equivalent to the left-hand side of Equation (\ref{poincare}) in the case
when $\hyp$ is of ideal type and the component $A$ 
is chosen to contain the dominant Weyl chamber.

One can speculate that there always exists a component $A$ such that 
$$P_{\hyp} (t) =  \prod_{i=1}^{n} (1 +t + t^2 + \dots + t^{m_i}).$$

One can also wonder if there exists a complex projective variety naturally associated
to $\hyp$ such that $P_{\hyp} (t^{\frac{1}{2}})$ is its Poincar{\' e} polynomial.
This would be the analogue of the regular nilpotent Hessenberg varieties.


\end{document}